\documentstyle[12pt]{article}

\begin{document}

\date{}

\title{\bf Fourier Knots}

\author{
Louis H. Kauffman \\
Department of Mathematics, Statistics and Computer
Science \\
University of Illinois at Chicago \\
851 South Morgan Street\\
Chicago, IL 60607-7045\\
}

 \maketitle

 \thispagestyle{empty}

 \section{Introduction}
This paper introduces the concept of Fourier knot. A Fourier knot is a knot
that is represented by a parametrized curve in three dimensional space such
that the three coordinate functions of the curve are each finite Fourier
series in the parameter.  That is, the knot can be regarded as the result
of independent vibrations in each of the coordinate directions with each of
these vibrations being a linear combination of a finite number of pure
frequencies.
\vspace{3mm}

The previously studied Lissajous knots \cite{1} constitute the case of a
single frequency in each coordinate direction. Not all knots are Lissajous
knots, and in fact the trefoil knot and the figure eight knot are the first
examples of non-Lissajous knots. The first section of this paper sketches
the proof that every tame knot is a Fourier knot. Subsequent sections give
robust examples of Fourier representations for the trefoil, the figure
eight and a class of knots that we call Fibonacci knots.  In the case of
the trefoil we have given a minimal Fourier representation in the sense
that it has single frequencies  in two of the coordinate directions and a
combination of frequencies in the third direction.  The paper ends by
pointing out the usual compact non-linear trigonometric formula for torus
knots, and raises the question of the finite Fourier series representations
for these knots.
\vspace{3mm}

On completing an early draft of this paper, we learned that an extensive
study of Fourier knots (there called Harmonic  Knots) has been carried out
by Aaron Trautwein in his 1995 PhD Thesis \cite{7} at the University of
Iowa under the direction of Jon Simon. While the independently obtained
results of the present paper are primarily illustrative of the idea of
Fourier knots, Trautwein's pioneering work establishes relationships
between the complexity of the harmonic representation and knot theoretic
indices such as crossing number and superbridge index.  The interested
reader should consult this work.
\vspace{3mm}

\noindent
{\bf Acknowledgment.}  The author thanks the National
Science Foundation  for  support of this research under
grant number  DMS-9205277.
\vspace{3mm}

\section{Every Knot is a Fourier Knot}
In considering problems about knots it is interesting to have an equation
for the knot or link under consideration.  By an equation for a knot I mean
the specification of  a parametrized curve in three dimensional space of
the form

				$$X(t) = A(t)$$
				$$Y(t) = B(t)$$
			         $$Z(t) =  C(t)$$

\noindent
where $A$, $B$, $C$ are smooth functions of the variable $t$ with (for a
knot) a total period of $P>0$, specifying an embedding of the circle into
three dimensional space.  Here the circle is the quotient space of the
interval from $0$ to $P$ , $[0,P]$, obtained by identifying the ends of the
interval to each other with the quotient topology.
\vspace{3mm}

A function $F(t)$ is said to be a finite Fourier series if it has the form

$$F(t) = A_{1}Cos(K_{1}T+L_{1}) + ...+A_{N}Cos(K_{N}T+L_{N})$$

\noindent
where $A_{1}$,$A_{2}$,...$A_{N}$ and $K_{1}$,$K_{2}$,...$K_{N}$ and
$L_{1}$,$L_{2}$,$L_{3}$,...$L_{N}$  are given constants and
$K_{1}$,$K_{2}$,...$K_{N}$ are each rational numbers. Of course, $F(t)$ can
also be expressed in terms of the $Sin$ function or as a combination of
$Sin$ and $Cos$ functions, since $Sin(X + \pi/2) = Cos(X). $
\vspace{3mm}

\noindent
{\bf Definition.}  A knot $K$ embedded in three dimensional space will be
called a {\em Fourier knot} if it has an equation (as described above) with
each of the functions $A$,$B$,$C$ a finite Fourier series.
\vspace{3mm}

\noindent
{\bf Definition.} A knot $K$ is said to be {\em tame} if every point $p$ of
$K$ in $R^{3}$ (Euclidean three-space) has a neighborhoodsuch that the
intersection of $K$ with that neighborhood is equivalent to a standard pair
of three-dimensional ball and diameter-arc of that ball.
\vspace{3mm}

It is well known that every tame knot can be represented topologically by
equations where $A$,$B$ and $C$ are smooth (i.e. infinitely differentiable)
functions. It then follows by standard approximation theorems for Fourier
series that $A$, $B$ and $C$ can be taken to be finite Fourier series.
Thus we have proved the
\vspace{2mm}

\noindent
{\bf Theorem.}  Every tame knot is (topologically equivalent to)  a
Fourier knot.
\vspace{3mm}

\section{Lissajous  Knots and the Arf Invariant}
One class of Fourier knots that have been studied are the Lissajous knots
\cite{1}.  In a Lissajous knot there is one term in each of the Fourier
series. Thus a Lissajous knot has the form

				$$X(t) = A_{1}Cos(K_{1}t+L_{1})$$
				$$Y(t) = A_{2}Cos(K_{2}t+L_{2})$$
				$$Z(t)  = A_{3}Cos(K_{3}t+L_{3})$$

In \cite{1} it is proved that the Arf invariant of a Lissajous knot is
neccessarily equal to zero.  This means that many knots are not Lissajous.
In particular the trefoil knot and the figure eight knot are not Lissajous
knots.  This leads to the question:  If a tame knot $K$ is not Lissajous,
what is the "simplest" representation of K in terms of finite Fourier
series?
\vspace{3mm}

In the next section we shall give a definite answer to this question in the
case of the trefoil knot, and a conjecture in the case of the figure eight
knot.  In all cases, when one answers this question there is also the
parallel question of obtaining Fourier equations for the knot K that are
robust in the sense that plots of these equations yield pleasing images
that can be explored geometrically.
\vspace{3mm}

Since knots of non-zero Arf invariant are neccessarily not Lissajous, it
will be useful for us to recall one definition of the Arf invariant. An
interested reader can apply this definition to find examples of Fourier
knots that are not Lissajous knots.  There are algebraic definitions of the
Arf invariant and a fundamental geometric definition as well. See \cite{2}
or \cite{3} for more details.
\vspace{3mm}

We recall an algebraic definition of the Arf invariant by first defining an
integer valued invariant, $a(K)$, associated to any oriented knot $K.$ The
invariant $a(K)$ is defined by the (recursive) equation (*)

$$ a(K_{+}) = a(K_{-}) = Lk(K_{0})$$

\noindent
where $K_{+}$, $K_{-}$ and $K_{0}$ are three diagrams that differ at a
single crossing as shown in Figure 1, and $Lk(K_{0})$ denotes the linking
number of the link of two components $K_{0}$.   $K_{+}$ and $K_{-}$  are
each knot diagrams, differing from each other by a single switched
crossing.  $K_{0}$ is obtained from either $K_{+}$ or $K_{-}$   by
smoothing that crossing. A {\em smoothing}  is accomplished by reconnecting
the strands at the crossing so that the arcs no longer cross over one
another (as shown in Figure 1).  Both the switching and the smoothing
operations can change the topological type of the diagram.  Smoothing
always replaces a knot by a link of two components. Thus $K_{0}$  is a such
a link.  By definition,  $a(K)$ is equal to zero if $K$ is topologically
equivalent to an unknotted circle.
\vspace{3mm}

\begin{figure}[htbp]
\vspace*{80mm}
\special{epsf=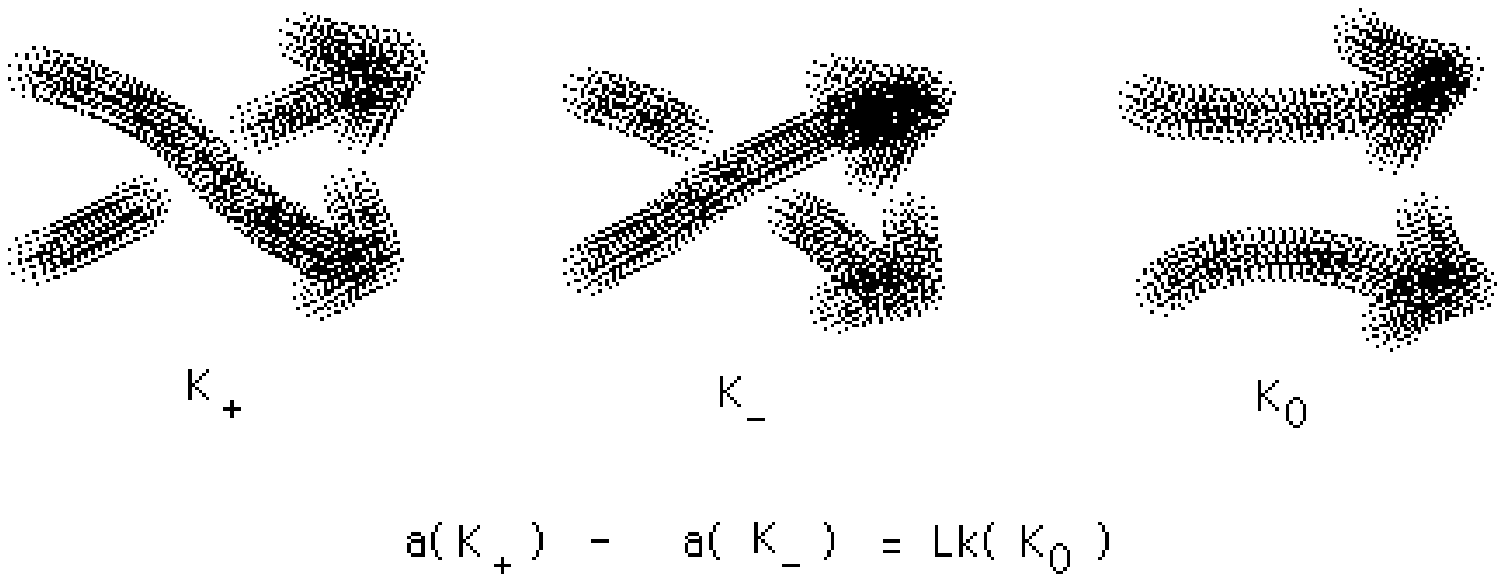}
\caption{Switching Relation}
\end{figure}

Note that in Figure 1 we have implicitly assigned signs of $+1$ or $-1$ to
the two types of oriented crossings. This number, $+1$ or $-1$, is called
the sign of the crossing. The linking number of a link L is defined by the
equation

$$Lk(L) = \Sigma_{p} e(p)/2$$

\noindent
where the summation runs over all the crossings in $K$ that are between two
components of $K.$ Crossings of any given component with itself are not
counted in this summation.
\vspace{3mm}

It is a (non-obvious) fact that the recursive equation (*) defines a
topological invariant of knots and links. It is, in fact part of a much
larger scheme of things.  For example it is the second coefficient of the
Conway polynomial.  One way to define the Arf invariant, $Arf(K)$, is by
the equation

$$Arf(K) = a(K) (mod 2).$$

Thus the Arf invariant of $K$ is either $0$ or $1$ depending upon the
parity of $a(K).$  See Figure 2 for a sample calculation of the Arf
invariant of the trefoil knot.
\vspace{3mm}

\begin{figure}[htbp]
\vspace*{80mm}
\special{epsf=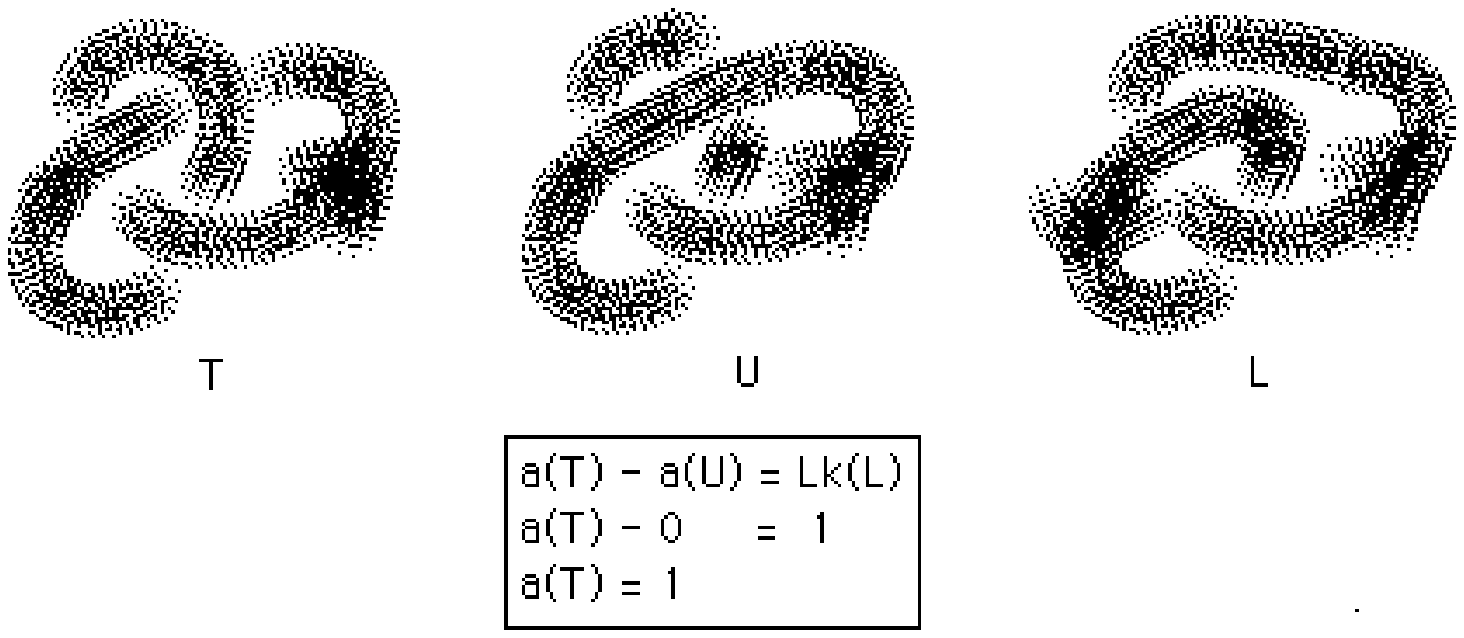}
\caption{a(Trefoil)}
\end{figure}

It is a remarkable fact that Lissajous knots have Arf invariant zero.
I do not know if every knot of vanishing Arf invariant is a Lissajous knot.
\vspace{3mm}

\section{A Fourier Trefoil Knot}
Consider the following equations

		$$x=Cos(2T),$$
		$$y=Cos(3T+ (1/2)),$$
		$$z=(1/2)Cos(5T+ (1/2)) + (1/2)Sin(3T + (1/2)).$$

\noindent
These equations define a trefoil knot, showing that the trefoil knot is a
Fourier knot where only one coordinate needs to be a combination of
frequencies. The proof that these equations give a trefoil knot is left to
the reader. One way to verify this is to use a computer to draw the
pictures in three dimensions and then examine the results. Figure 3
illustrates a computer drawing of this Fourier trefoil.  The drawing
illustrates what I mean by a robust representation of the knot. The knot
does not come ambiguously close to itself, and the form of the drawing is
aesthetically pleasing.
\vspace{3mm}

\begin{figure}[htbp]
\vspace{80mm}
\special{epsf=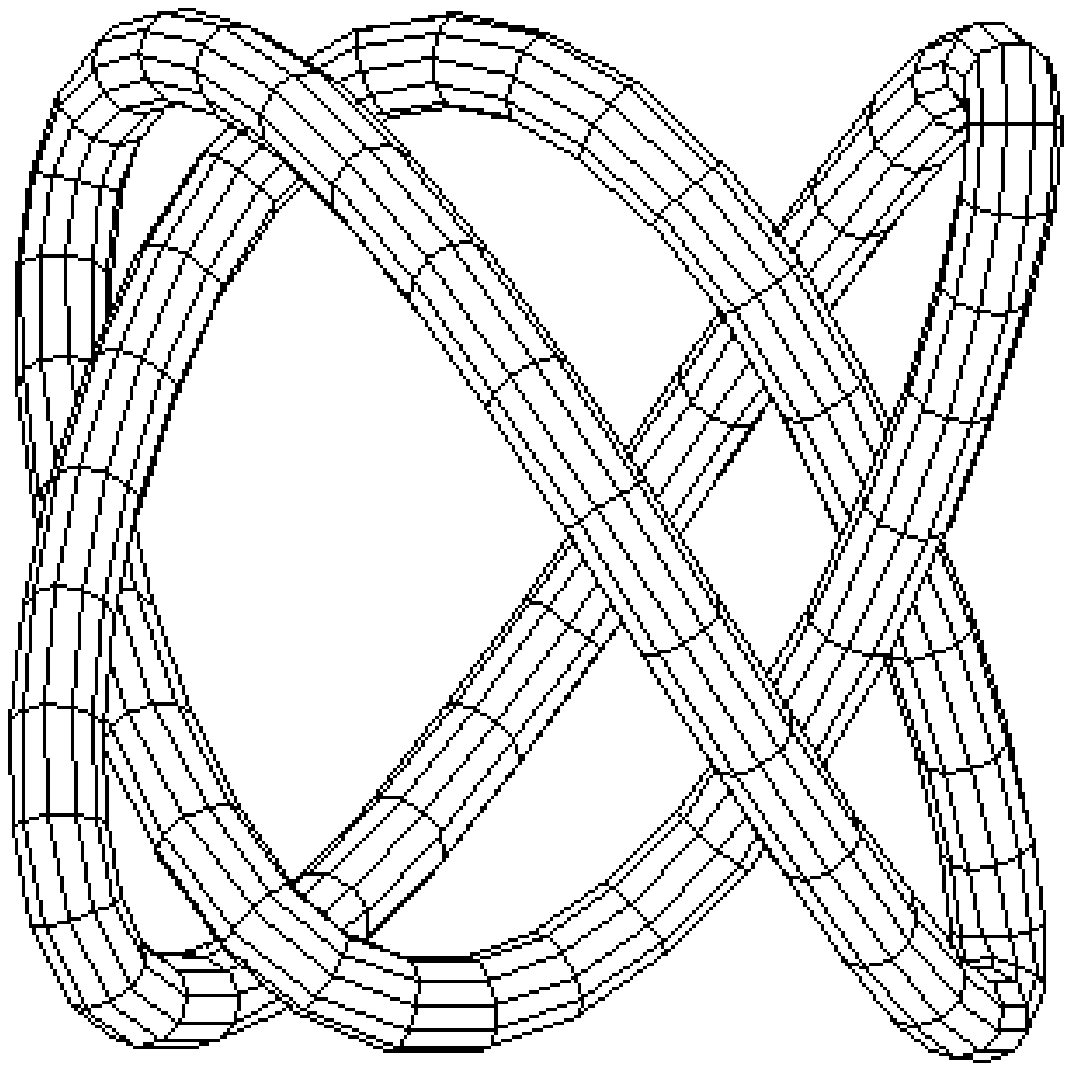}
\caption{The Fourier Trefoil}
\end{figure}

The author wishes to acknowledge Lynnclaire Dennis \cite{4} for inspiring
him to search for the Fourier trefoil. In her book "The Pattern" Ms. Dennis
draws a picture of a knot (the Pattern knot)  that closely resembles our
Fourier trefoil. In projection the Pattern knot looks like a Lissajous
figure with frequencies $2$ and $3$ and the Pattern knot is a trefoil knot.
This led to trying combinations of frequencies for the third coordinate,
and eventually to the equations above with pure frequencies $2$ and $3$ in
two directions and the combination of frequencies $5$ and $3$ in the third
direction.  The Pattern knot is more spherically symmetrical than the
Fourier trefoil, and does not have an obvious equation.
\vspace{3mm}

I would also like to mention an experiment that I performed with the
Fourier trefoil in the form of a (hand-drawn)  diagram corresponding to the
knot in Figure 3. I gave this diagram as input to Ming, a knot energy
program written by Ying-Quing Wu at the University of Iowa.  (A
diagrammatic interface for Ming was written by Milana Huang at the
Electronic Visualisation Lab at the University of Illinois).  Ming sets the
knot on a descending energy trajectory, following Jon Simon's energy
\cite{5}  for the knot. The result of this experiment is that the flat knot
diagram quickly unfurls into a three dimensional geometry very similar to
the Fourier trefoil and nearly stabilizes in this form. Then slowly the
knot moves off this slightly higher energy level and settles into the
familiar symmetry of the (empirically) known energy minimum for the trefoil
knot.  Thus there appears to be a "point of inflection" in this particular
way of descending to minimum energy for the trefoil knot. This experiment
points to a wide range of possible explorations, investigating the gradient
descent for knot energy from particular starting configurations for a knot.
\vspace{10mm}

\section{A Fourier Figure Eight Knot}
The following equations describe a figure eight knot.

			$$x=Cos(t) + Cos(3t),$$
			$$y=.6Sin(t) + Sin(3t),$$
			$$z=.4Sin(3t)  - Sin(6t).$$
\noindent
See Figure 4.

\begin{figure}[htbp]
\vspace*{100mm}
\special{epsf=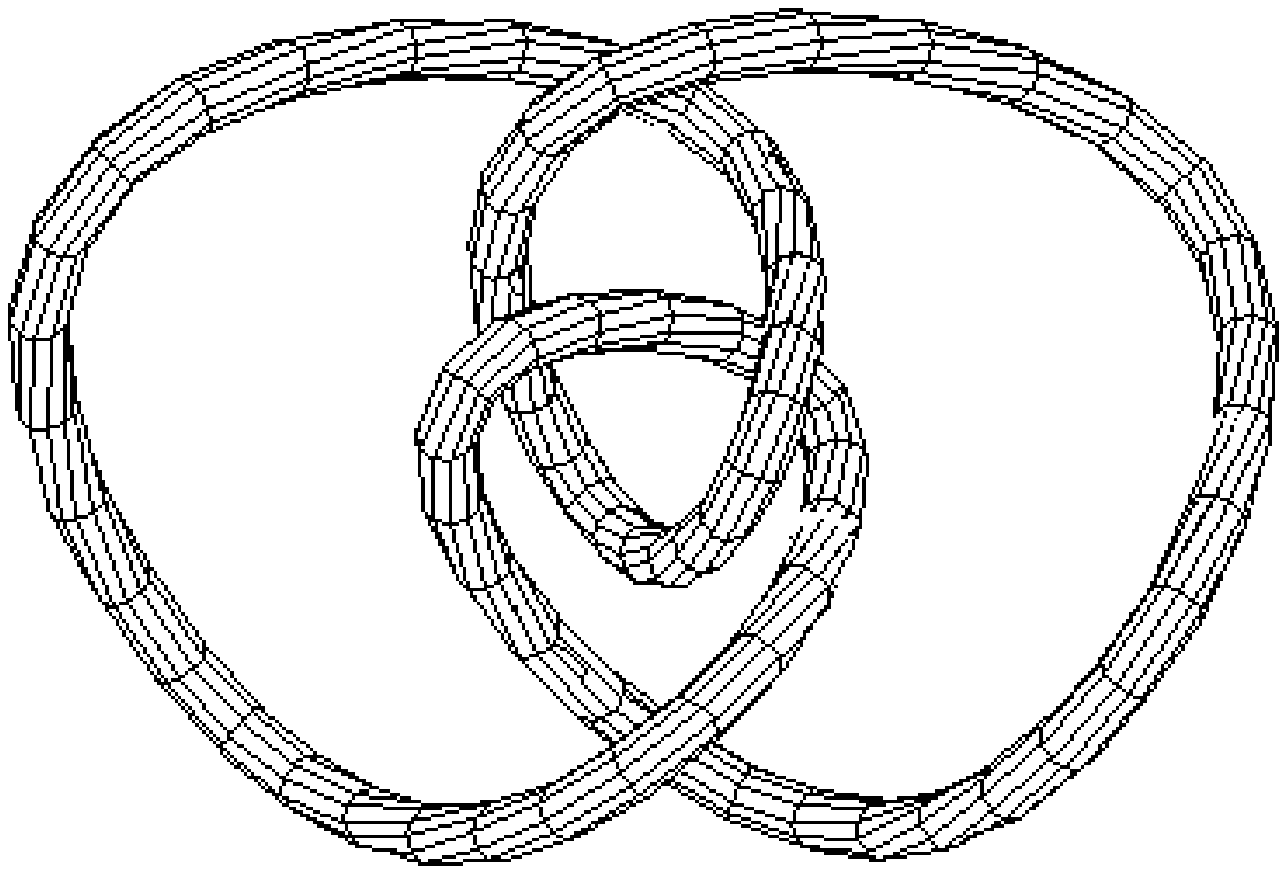}
\caption{The Fourier Figure Eight Knot}
\end{figure}

\noindent
I do not know if there is a simpler Fourier representation for this knot.
\vspace{3mm}

\section{\bf A Series of Fibonacci Fourier Knots}
Recall the Fibonacci series

$$1,1,2,3,5,8,13,21,34,55,89,144,....$$

The $n$-th term, $f_{n}$,  of the series is equal to the sum of the
previous two terms, and $f_{1} = f_{2}= 1.$ Consider the equations

			$$x=Cos(8T),$$
			$$y=Cos(13T+.5),$$
			$$z=.5Cos(21T+.5)+.5Sin(13T +.5),$$

\noindent
and more generally

			$$x=Cos(f_{n}T),$$
			$$y=Cos(f_{n+1}T+.5),$$
			$$z=.5Cos(f_{n+2}T+.5)+.5Sin(f_{n+1}T +.5).$$

\noindent
The last set of equations defines a knot that we shall dub $F(n)$, the
$n$-th Fibonacci knot.  Thus $F(3)$ is the Fourier trefoil of Section 3,
and the first equations we have written in this section denote $F(6)$.  In
Figure 5 we illustrate a computer drawing of $F(6).$
\vspace{3mm}

\begin{figure}[htbp]
\vspace*{130mm}
\special{epsf=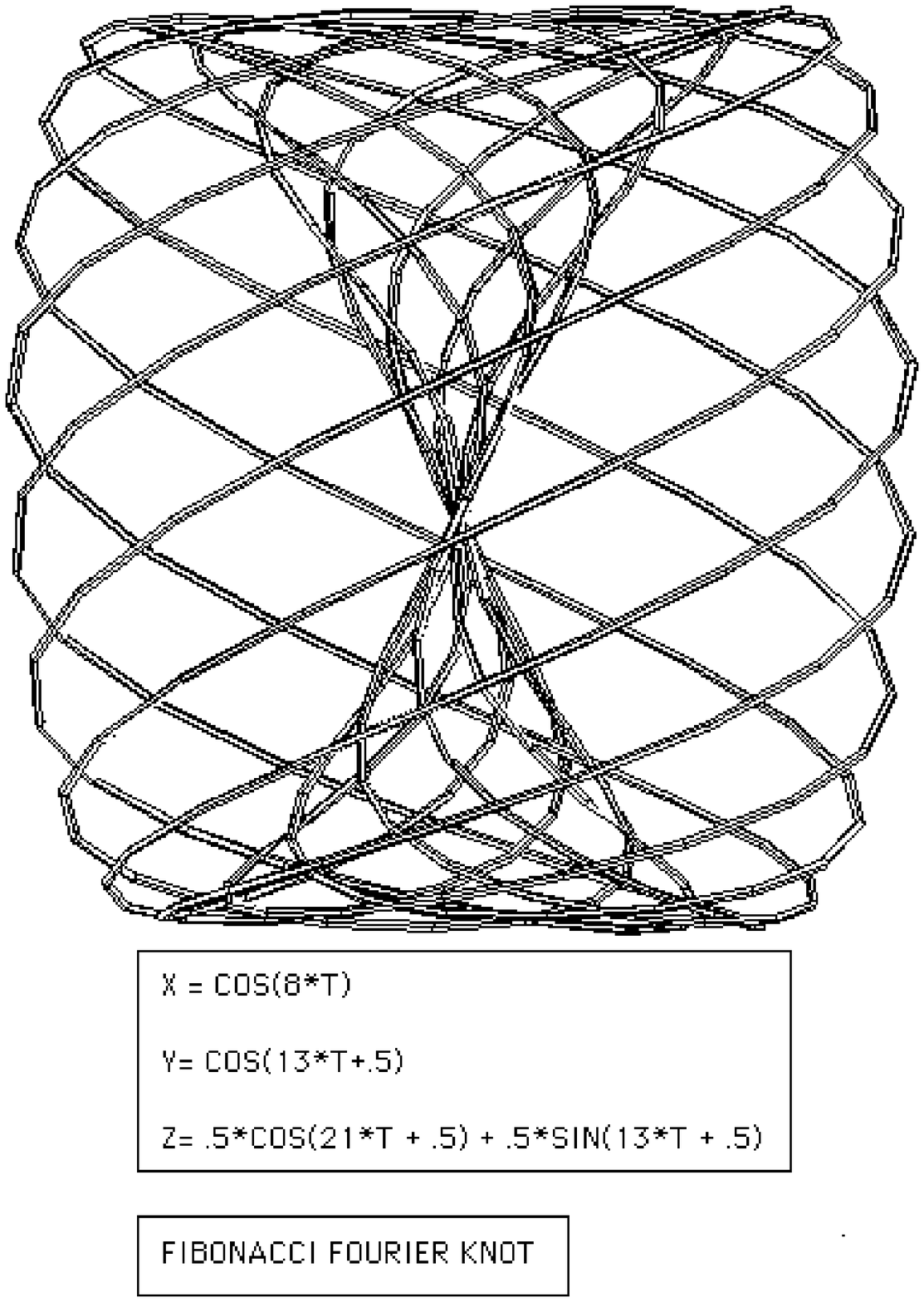}
\caption{A Fibonacci Fourier Knot}
\end{figure}

The Fibonacci Fourier knots provide a strong class of knots for
investigation using computer graphics.
\vspace{3mm}

\section{Torus Knots}
Recall that a knot that winds $P$ times around a torus in one direction and
$Q$ times in the other direction -  a torus knot of type $(P,Q)$ - has the
equation

			$$x=Cos(T)(1 + .5Cos((Q/P)T)),$$
			$$y= Sin(T)(1 +.5Cos((Q/P)T)),$$
			$$z= .5Sin((Q/P)T).$$

\noindent
 Now use the trigonometric identities
$$cos(a)cos(b) = .5(cos(a+b) +cos(a-b)),$$
$$sin(a)cos(b) = .5(sin(a+b) + sin(a-b)).$$
The equations above then become

$$x=cos(T) + .25cos((1+Q/P)T) + .25cos((1-Q/P)T)$$
$$y=sin(T) + .25*sin((1+Q/P)T) + .25*sin((1-Q/P)T)$$
$$z= .5sin((Q/P)T)$$

\vspace{3mm}

\noindent
Thus torus knots are Fourier knots, and we can ask if these are simplest
Fourier representations for torus knots.
The parametrization shown above appears in \cite{6}. I am indebted to Peter
Roegen for pointing this out.

 \end{document}